\documentclass[11pt]{amsart}

\usepackage{amssymb}

\evensidemargin\oddsidemargin

\begin{document}
\baselineskip=18pt
\setcounter{page}{1}
    
\newtheorem{Conj}{Conjecture\!\!}
\newtheorem{Theo}{Theorem\!\!}
\newtheorem{Lemm}{Lemma}
\newtheorem{Rem}{Remark\!\!}

\renewcommand{\theConj}{}
\renewcommand{\theRem}{}
\renewcommand{\theTheo}{}

\def\a{\alpha}
\def\b{\beta}
\def\B{{\bf B}} 
\def\C{{\bf C}} 
\def\GG{{\mathcal{G}}} 
\def\KK{{\mathcal{K}}} 
\def\LL{{\mathcal{L}}} 
\def\SS{{\mathcal{S}}}
\def\UU{{\mathcal{U}}}
\def\ca{c_{\a}}
\def\ka{\kappa_{\a}}
\def\coa{c_{\a, 0}}
\def\cua{c_{\a, u}}
\def\cL{{\mathcal{L}}} 
\def\Ea{E_\a}
\def\eps{{\varepsilon}} 
\def\esp{{\mathbb{E}}} 
\def\Ga{{\Gamma}} 
\def\G{{\bf \Gamma}} 
\def\e{{\rm e}}
\def\i{{\rm i}}
\def\L{{\bf L}}
\def\lbd{\lambda}
\def\lcr{\left[}
\def\lpa{\left(}
\def\lva{\left|}
\def\M{{\bf M}}
\def\NN{{\mathbb{N}}} 
\def\pb{{\mathbb{P}}}
\def\rl{{\mathbb{R}}}
\def\rpa{\right)}
\def\rcr{\right]}
\def\rva{\right|}
\def\W{{\bf W}}
\def\X{{\bf X}}
\def\XX{{\mathcal X}}
\def\YY{{\mathcal Y}}
\def\U{{\bf U}}
\def\V{{\bf V}_\a}
\def\Un{{\bf 1}}
\def\Z{{\bf Z}}
\def\A{{\bf A}}
\def\AA{{\mathcal A}}
\def\hAA{{\hat \AA}}
\def\hL{{\hat L}}
\def\hT{{\hat T}}

\def\claw{\stackrel{d}{\longrightarrow}}
\def\elaw{\stackrel{d}{=}}
\def\qed{\hfill$\square$}

\title[Proof of Bondesson's conjecture]{A proof of Bondesson's conjecture on stable densities}

\author[Pierre Bosch]{Pierre Bosch}

\address{Laboratoire Paul Painlev\'e, Universit\'e Lille 1, Cit\'e Scientifique, F-59655 Villeneuve d'Ascq Cedex. {\em Email}: {\tt pierre.bosch@math.univ-lille1.fr}}

\author[Thomas Simon]{Thomas Simon}

\address{Laboratoire Paul Painlev\'e, Universit\'e Lille 1, Cit\'e Scientifique, F-59655 Villeneuve d'Ascq Cedex. Laboratoire de Physique Th\'eorique et Mod\`eles Statistiques, Universit\'e  Paris-Sud, F-91405 Orsay. {\em Email}: {\tt simon@math.univ-lille1.fr}}



\begin{abstract} We show that positive $\a-$stable densities are hyperbolically completely monotone if and only if $\a\le 1/2.$ This gives a positive answer to a question raised by L.~Bondesson in 1977.
\end{abstract}

\maketitle

\section{Introduction and statement of the result}

This paper is a contribution to the analytical structure of positive stable densities.  Up to multiplicative normalization, the latter are defined through their Laplace transform 
$$\int_0^\infty e^{-\lbd x} f_\a(x)\, dx\; = \; e^{-\lbd^\a},$$
where $\a\in (0,1)$ is a self-similarity parameter. We will denote by $\Z_\a$ the underlying positive $\a-$stable random variable. The density $f_\a$ is not explicit except in the case $\a = 1/2,$ where it writes 
\begin{equation}
\label{Expl}
f_{1/2} (x)\; =\; \frac{1}{2\sqrt{\pi x^3}}\, e^{-\frac{1}{4x}}.
\end{equation}
When $\a$ is rational, various explicit factorizations of the random variable $\Z_\a$ are available in terms of Beta and Gamma random variables - see Section 2 in \cite{TS2} and the references therein. The importance of the densities $f_\a$ stems from the fact that they serve as building blocks in order to construct all strictly stable densities on the line - see Section 7 in \cite{TS2}. Since their introduction by L\'evy in the late twenties, real stable densities are of constant use in probability and statistics. We refer to the classic monograph \cite{Z} for a presentation of such densities.

In this paper, we are interested in the following property. A density function $f$ on $\rl^+$ is called hyperbolically completely monotone (HCM for short) if it is the pointwise limit of functions of the type
$$x\;\mapsto\;C\, x^{\b - 1}\,\prod_{i = 1}^N\, (1+ xy_i)^{-\gamma_i},$$
where all parameters are positive. The denomination comes from the equivalent definition that for every $u > 0,$ the function $ f(uv)f(u/v)$ is completely monotone in the variable $w=v + 1/v.$ This notion was introduced by Thorin and Bondesson in the late seventies in order to show the infinite divisibility of densities on the half-line having untractable Laplace transform, like the log-normal or the generalized Pareto. It was shown that HCM densities are all additive convolutions of Gamma densities with different parameters, a property called generalized Gamma convolution (GGC for short), which entails infinite divisibility. Another important fact is the stability of the HCM property by power transformations of absolute value greater than one. We refer to the monograph \cite{B} for the classic account on these two notions. See also \cite{BJSP} for a short presentation and an extension to the line. The GGC property of a positive random variable amounts to the fact that its log-Laplace exponent is a Pick-Nevanlinna function, and with this characterization it is easy to see - see Example 3.2.1 in \cite{B} - that positive stable laws are all GGC's. On the other hand, the stringent HCM property is more difficult to study for positive stable densities, because it is read off on the density itself, which is here not explicit in general.
\begin{Theo} The density $f_\a$ is {\em HCM} if and only if $\a \le 1/2.$
\end{Theo}
This result was conjectured in 1977 by Bondesson. We refer to  \cite{BZW} pp. 54-55, Sections 5.7, 5.8 and 7.2 in \cite{B}, Remark 3 in \cite{BJSP} and to the whole paper \cite{B99} for several reasons, partly numerical, supporting this hypothesis. Notice that the only if part is easy to establish because the random variable $\Z_\a^{-1}$ is not infinitely divisible, and hence cannot have a HCM density, when $\a > 1/2.$ The if part had been observed in Example 5.6.2 of \cite{B} for the case $\a = 1/n$ with $n\ge 2,$ thanks to a factorization of $\Z_{1/n}^{-1}$ in terms of Gamma random variables - see \cite{W} - which will be discussed in Section 3 below.  Unfortunately, such a factorization cannot hold when $\a$ is no more the reciprocal of an integer - see Example 7.2.5 in \cite{B}. In \cite{TS1}, it was shown with the help of Beta-Gamma products that $f_\a$ is hyperbolically monotone (that is, the
function $ f_\a(uv)f_\a(u/v)$ is non-increasing in the variable $w=v + 1/v$) if and only if $\a \le 1/2.$ Other partial results supporting Bondesson's conjecture were obtained in \cite{JS}, including the facts that the independent products $\Z_\a \times \Z_\a^{-1}$ and $\G_1^{1/\a}\times \Z_\a$ have HCM densities if and only if $\a \le 1/2.$ Finally, a proof of this conjecture in the case $\a\not\in (1/4, 1/3)$ was recently announced in \cite{Fouras}.

Our method to solve this problem relies on a certain representation of $\Z_\a^{-1}$ as an infinite product of independent renormalized Beta random variables, which was recently observed in \cite{LS} in the framework of more general stable functionals, and is a simple consequence of Malmsten's formula for the Gamma function. In the case $\a\le 1/2,$ this representation entails that $\Z_\a^{-1}$ factorizes into the product of a Gamma random variable $\G_\a$ and another infinite Beta product. Using a recent key-result by Bondesson \cite{BJTP} on the stability of the GGC property by independent multiplication of random variables, we can then deduce, without much effort, that this latter factorization has indeed an HCM density. 

\section{Proof of the Theorem}

Throughout this note, we will denote by $\B_{a,b}$ and $\G_c$ the usual Beta and Gamma random variables, with respective densities
$$f_{\B_{a,b}}^{}(x)\; =\; \frac{\Ga(a+b)}{\Ga(a)\Ga(b)} \, x^{a-1}(1-x)^{b-1} \, \Un_{(0,1)}(x)\qquad\mbox{and}\qquad f_{\G_c}^{}(x)\; =\; \frac{1}{\Ga(c)} \, x^{c-1}\, e^{-x} \, \Un_{(0,+\infty)}(x).$$
It will be implicitly assumed that all quotients and products of given random variables are independent. The following key-lemma can be viewed as a generalization of Lemma 1 in \cite{TS1}.

\begin{Lemm} The density of the product
$$\G_c\;\times\; \B_{a_1, b_1}\;\times\;\cdots\;\times\; \B_{a_n, b_n}$$
is {\em HCM} for every $n\ge 1, a_i, b_i > 0$ and $c < \min\{a_i\}.$
\end{Lemm}

\proof Set $f,g$ for the respective densities of $\B_{a_1, b_1}\times\cdots\times \B_{a_n, b_n}$ and $\G_c\times \B_{a_1, b_1}\times\cdots\times \B_{a_n, b_n}.$ The multiplicative convolution formula and the change of variable $y=(z+1)^{-1}$ show that
$$g(x) \; = \; \frac{x^{c-1}}{\Ga(c)}\, \int_0^1 e^{-xy^{-1}}\, y^{-c} f(y)\, dy\; =\; \frac{x^{c-1}e^{-x}}{\Ga(c)}\, \int_0^\infty e^{-xz}\, (z+1)^{c-2} f((z+1)^{-1})\, dz.$$
On the other hand, the condition $c < \min\{a_i\}$ entails easily from a comparison of  the Mellin transforms that $y\mapsto y^{-c} f(y)$ is, up to normalization, the density of $\B_{a_1-c, b_1}\times\cdots\times \B_{a_n-c, b_n}.$ Hence, the function $z\mapsto 
(z+1)^{c-2} f((z+1)^{-1})$ is up to normalization the density of 
$$\B_{a_1-c, b_1}^{-1}\times\cdots\times \B_{a_n-c, b_n}^{-1}\; -\; 1.$$
The density of $\B_{a_i-c, b_i}^{-1} - 1$ is 
$$\frac{\Ga(a_i+b_i -c)}{\Ga(a_i -c)\Ga(b_i)} \,\times\, \frac{x^{b_i-1}}{(1+x)^{a_i+b_i-c}}$$
and is HCM, so that the law of $\B_{a_i-c, b_i}^{-1}$ is a GGC for every $i=1\ldots n.$ By the main result of \cite{BJTP}, the law of $\B_{a_1-c, b_1}^{-1}\times\cdots\times \B_{a_n-c, b_n}^{-1} - 1$ is hence a GGC as well. A combination of Theorem 5.4.1 and Property (ii) p.68 in \cite{B} shows finally that $g$ is HCM.
 
\endproof

\begin{Rem} {\em By stability of the HCM property under weak convergence - see Theorem 5.1.3 in \cite{B}, the condition on $c$ can be relaxed into $c\le\min\{a_i\}.$ A perusal of the proof of Lemma 1 in \cite{TS1} shows also that the product $\G_c\times\B_{a,b}$ has a HCM density as soon as $a+b\ge c.$ We believe that this latter condition is necessary.}
\end{Rem}

The next lemmas provide infinite Beta product representations for positive stable and Gamma random variables, which are interesting in their own right. The proof is essentially the same as that of the main result in \cite{LS}, but we provide all details for the sake of completeness.

\begin{Lemm} For every $\a\in (0,1),$ one has the a.s. convergent factorization
$$\Z_\a^{-1}\; \elaw\; e^{\gamma(1-\a^{-1})}\,\times\,\prod_{n=0}^\infty a_{n}\, \B_{\a+n\a,1-\a}$$
where $\gamma$ is Euler's constant, $\psi$ is the digamma function and $a_n = e^{\psi(1+n\a) - \psi(\a +n\a)}.$
\end{Lemm}

\proof The proof is based on the well-known representation
\begin{equation}
\label{Malm}
\frac{\Ga (a+s)}{\Ga(a)}\; =\; \exp\lcr \psi(a) s\; +\; \int_{-\infty}^{0} (e^{sx} - 1 -sx)\, \frac{e^{-a\vert x\vert}}{\vert x\vert (1-e^{-\vert x\vert})}\, dx\rcr
\end{equation}
for every $a,s > 0,$ which is easily obtained from Malmsten's formula for the Gamma function - see e.g. 1.9(1) p.21 in \cite{EMOT}. Making some simplifications, we deduce
\begin{eqnarray*}
\esp[\Z_\a^{-s}] & = & \frac{\Ga (1+s\a^{-1})}{\Ga(1+s)}\\
& =&  \exp\lcr \gamma(1-\a^{-1})s\; +\; \int_{-\infty}^{0} (e^{sx} - 1 -sx)\, \frac{e^{-\a\vert x\vert} (1-e^{-(1-\a)\vert x\vert})}{\vert x\vert (1-e^{-\vert x\vert}) (1-e^{-\a\vert x\vert})}\, dx\rcr \\
& = & \exp\lcr \gamma(1-\a^{-1})s\; +\;\sum_{n=0}^\infty \int_{-\infty}^{0} (e^{sx} - 1 -sx)\, \frac{e^{-(\a+n\a)\vert x\vert} (1-e^{-(1-\a)\vert x\vert})}{\vert x\vert (1-e^{-\vert x\vert})}\, dx\rcr,
\end{eqnarray*}
where the third equality follows from Fubini's theorem. On the other hand, it is well-known and easy to see from (\ref{Malm}) that
$$\esp[\B_{\a+n\a,1-\a}^{s}]\; =\; \frac{\Ga (\a+n\a +s)\Ga(1+n\a)}{\Ga(1+n\a+s)\Ga(\a+n\a)}\; =\;  \exp\lcr \int_{-\infty}^{0} (e^{sx} - 1)\, \frac{e^{-(\a+n\a)\vert x\vert} (1-e^{-(1-\a)\vert x\vert})}{\vert x\vert (1-e^{-\vert x\vert})}\, dx\rcr,$$
which also entails
$$\esp[\log \B_{\a+n\a,1-\a}]\; =\;-\int_{-\infty}^{0} \frac{e^{-(\a+n\a)\vert x\vert} (1-e^{-(1-\a)\vert x\vert})}{(1-e^{-\vert x\vert})}\, dx \; =\; \psi(\a+n\a)-\psi(1+n\a)$$
(see e.g. Formula 1.7.2(14) in \cite{EMOT} for the second equality). Finally, we observe that the martingale
$$X_n\;=\;\sum_{i=0}^n \lpa \log \B_{\a+i\a,1-\a} - \esp[\log \B_{\a+i\a,1-\a}]\rpa$$
where all summands are assumed independent, converges a.s. by the martingale convergence theorem, since the variances are uniformly bounded by
\begin{eqnarray*}
\int_0^{\infty} \frac{x\,e^{-\a x} (1-e^{-(1-\a)x})}{(1-e^{-x})(1-e^{-\a x})}\, dx & < & +\infty.\end{eqnarray*}

\noindent
Putting everything together and identifying the Mellin transforms shows the required factorization
$$\Z_\a^{-1}\; \elaw\; e^{\gamma(1-\a^{-1})} \lpa \lim_{n\to +\infty} e^{X_n}\rpa\; \elaw\; e^{\gamma(1-\a^{-1})}\,\times\,\prod_{n=0}^\infty a_{n}\, \B_{\a+n\a,1-\a}.$$

\endproof

\begin{Lemm} For every $a, b > 0,$ one has the a.s. convergent factorization
$$\G_a\; \elaw\; e^{\psi(a)}\,\times\,\prod_{n=0}^\infty b_{n}\, \B_{a+nb,b}$$
with $b_n = e^{\psi(a +b+nb) - \psi(a +nb)}.$
\end{Lemm}

\proof The proof is an analogous consequence of (\ref{Malm}), which can be rewritten
$$\esp[\G_a^s]\; =\; \frac{\Ga (a+s)}{\Ga(a)}\; =\; \exp\lcr \psi(a) s\; +\; \int_{-\infty}^{0} (e^{sx} - 1 -sx)\, \frac{e^{-a\vert x\vert}(1-e^{-b\vert x\vert})}{\vert x\vert (1-e^{-\vert x\vert})(1-e^{-b\vert x\vert})}\, dx\rcr$$
for every $a,b >0.$ We omit the details. 

\endproof

We can now proceed to the proof of the Theorem. Recalling the discussion made in the introduction, we need to show that $\Z_\a^{-1}$ has a HCM density for every $\a < 1/2.$ Applying the elementary factorization
\begin{equation}
\label{Elem}
\B_{a,b+c}\; \elaw\;\B_{a,b}\,\times\,\B_{a+b,c},
\end{equation}
we see from Lemmas 2 and 3 that
\begin{eqnarray*}
\Z_\a^{-1} & \elaw & e^{\gamma(1-\frac{1}{\a})}\,\times\,\prod_{n=0}^\infty e^{\psi(1+n\a) - \psi(\a +n\a)}\, (\B_{\a+n\a,\a}\times \B_{2\a+n\a,1-2\a})\\
& \elaw & e^{\gamma(1-\frac{1}{\a})}\,\times\lpa\prod_{n=0}^\infty e^{\psi(2\a+n\a) - \psi(\a +n\a)} \,\B_{\a+n\a,\a}\rpa\times\lpa \prod_{n=0}^\infty e^{\psi(1+n\a) - \psi(2\a +n\a)}\,\B_{2\a+n\a,1-2\a}\rpa\\
& \elaw & e^{\gamma(1-\frac{1}{\a}) -\psi(\a)}\,\times\, \G_\a \,\times\lpa \prod_{n=0}^\infty e^{\psi(1+n\a) - \psi(2\a +n\a)}\,\B_{2\a+n\a,1-2\a}\rpa
\end{eqnarray*}
with again an a.s. convergent infinite product on the right-hand side, by the martingale convergence theorem. This shows that $\Z_\a^{-1}$ is the a.s. limit of
$$\KK_n\,\times\, \G_\a \;\times\;\B_{2\a,1-2\a}\;\times\;\cdots\;\times\;\B_{2\a+n\a,1-2\a}$$
as $n\to +\infty,$ for some deterministic renormalizing constant $\KK_n.$ Applying  Lemma 1 and Theorem 5.1.3 in \cite{B} concludes the proof.

\qed

\section{Further remarks}

\subsection{Connection with Gamma product representations} The infinite Beta product representation of $\Z_\a^{-1}$ obtained in Lemma 2 can be viewed as an extension of the factorization
$$\Z_{\frac{1}{n}}^{-1}\; \elaw\; n^n \times\, \G_{\frac{1}{n}}\,\times\,\cdots\,\times\, \G_{\frac{n-1}{n}},$$
which is due to Williams \cite{W}. A combination of Lemmas 2 and 3 shows for example that 
$$\Z_{\frac{1}{2}}^{-1}\; \elaw\; e^{-\gamma}\,\times\,\prod_{n=0}^\infty a_{n}\, \B_{\frac{n+1}{2},\frac{1}{2}}\; \elaw\; e^{-\gamma-\psi(\frac{1}{2})}\,\times\,\G_{\frac{1}{2}}\; =\; 4 \,\G_{\frac{1}{2}}$$
where the third equality comes from Formula 1.7.1(12) in \cite{EMOT}. This, of course, matches also the explicit expression (\ref{Expl}) for the density of $\Z_{\frac{1}{2}}$. Similarly, Lemma 2 and a repeated use of (\ref{Elem}) and Lemma 3 shows that
$$\Z_{\frac{1}{n}}^{-1}\; \elaw\; e^{-(n-1)\gamma-\psi(\frac{1}{n}) - \cdots-\psi(\frac{n-1}{n})} \,\times\, \G_{\frac{1}{n}}\,\times\,\cdots\,\times\, \G_{\frac{n-1}{n}}\; \elaw\; n^n \times\, \G_{\frac{1}{n}}\,\times\,\cdots\,\times\, \G_{\frac{n-1}{n}}.$$
As mentioned in the introduction, the representation of $\Z_\a^{-1}$ as a Gamma product (finite or infinite) is impossible when $\a$ is not the reciprocal of an integer, since then the law of $\log \Z_\a$ is not an extended GGC anymore - see Section 7.2 in \cite{B} for details. Other factorizations involving Beta and Gamma random variables are discussed in Section 2 of \cite{TS2} when $\a$ is rational, involving powers. The latter are not quite adapted to the HCM problem, which heavily depends on powers - see next paragraph. The fact that an appropriate extension of Williams' factorization should be searched for in order to solve the HCM problem for positive stable densities was suggested in the introduction of \cite{B99}. 

\subsection{Extension to power transformations} Our result implies that the density of $\Z_\a^q$ is HCM as soon as $\a\le 1/2$ and $\vert q\vert \ge 1$ - see Chapter 5 in \cite{B} p. 69. Observe on the other hand that any power transformation $\Z_\a^q$ does not have an HCM density when $\a > 1/2,$ since the latter density is not even hyperbolically monotone - see the main result of \cite{TS1}. It is a natural question to ask whether the optimal power exponent can be lowered in the case $\a\le 1/2$ and in this respect, the following was formulated in \cite{Bo}.

\begin{Conj} The density of $\Z_\a^q$ is {\em HCM} if and only if $\a\le 1/2$ and $\vert q\vert \ge \a/(1-\a).$
\end{Conj}

The only if part of this conjecture is also established in \cite{Bo}, together with the related fact that the density of the quotient $(\Z_\a \times \Z_\a^{-1})^q$ is HCM when $\a\le 1/2$ and $\vert q\vert \ge \a/(1-\a).$ The further factorization
$$\Z_\a^{-\a}\; \elaw\; e^{\gamma(\a-1)}\,\times\,\prod_{n=0}^\infty b_{n}\, \B_{1+\frac{n}{\a},\frac{1}{\a}-1}$$
where $b_n = e^{\psi(1+\frac{n}{\a}) - \psi(\frac{n+1}{\a})},$ obtained similarly as in Lemma 2, can be used to deduce from Theorem 2 in \cite{BS} and the main result of \cite{BJTP} that the law of $\Z_\a^q$ is a GGC for every $\a\in (0,1)$ and $q\ge \a/2.$ This factorization seems however helpless to tackle the HCM problem for lower powers of $\Z_\a$, which we postpone to future research.

%\bigskip
  
%\noindent
%{\bf Acknowledgement.} We thank Lennart Bondesson for having suggested the problem.

\end{document}